\magnification\magstep1 \baselineskip = 18pt \def\n{\noindent}
\def\vp{\varepsilon}

\magnification\magstep1 \baselineskip = 18pt \def\n{\noindent}

\def\pf{\medskip{\noindent{\bf Proof. }}} \def \rat{ {\rm Q}\kern-.65em
{}^{{}_/ }}

 \baselineskip = 18pt \def\n{\noindent} \overfullrule = 0pt
 \def\qed{{\hfill{\vrule height7pt width7pt depth0pt}\par\bigskip}}

 \magnification\magstep1

\def \rat{ {\rm Q}\kern-.65em {}^{{}_/ }}

\overfullrule = 0pt \def\pf{\medskip{\noindent{\bf Proof. }}}

\overfullrule = 0pt \def\pf{\medskip{\noindent{\bf Proof. }}}

\centerline{\bf Projections from a von~Neumann algebra} \centerline{\bf
onto a subalgebra} \bigskip \centerline{by}\bigskip \centerline{Gilles
Pisier\footnote*{Supported in part by the NSF}} \centerline{Texas A\&M
University} \centerline{and} \centerline{Universit\'e Paris VI}\bigskip

\n {\bf Introduction.} This paper is mainly devoted to the following
question:\ Let $M,N$ be von~Neumann algebras with $M\subset N$, if
there is a completely bounded ($c.b.$ in short) projection $P\colon
\ N\to M$, is there automatically a contractive projection $\widetilde
P\colon \ N\to M$?

We give an affirmative answer with the only restriction that $M$ is
assumed semi-finite. At the time of this writing, the case when the
subalgebra $M$ is a type~III factor seems unclear, although this might
be not too hard to deduce from our results using crossed product
techniques from the Tomita-Takesaki theory with which we are not
familiar.

If $N=B(H)$, a positive answer (without any restriction on $M$) was
given in [P1, P2] (and independently in [CS]). I am grateful to
Eberhard Kirchberg for mentioning to me that a more general statement
might be true. It should be mentioned that the above question seems
open if ``completely bounded'' is replaced by ``bounded'' in the
assumption on the projection $P$. For more results in this direction,
see [P3] and [HP2]. We should recall that, by a classical result of
Tomiyama [T],   every norm one projection $P$ from $N$ onto $M$
necessarily is a conditional expectation and in particular is
completely positive.  In the second part of the paper we give an
interpolation theorem which generalizes a result in [P1],  as follows.
Let $N$ be a von~Neumann algebra equipped with a normal semi-finite
faithful trace $\varphi$. Let us denote by $L_p(\varphi)$ the
noncommutative $L_p$-space associated to $(N,\varphi)$ in the usual
way. Fix $n\ge 1$. Let us denote by $A_0$ (resp.  $A_1$) the space
$N^n$ equipped with the norms $$\eqalign{\|(x_1,\ldots, x_n)\|_{A_0} &=
\left\|\left(\sum x_ix^*_i\right)^{1/2}\right\|_N\cr \|(x_1,\ldots,
x_n)\|_{A_1} &= \left\|\left(\sum x^*_i x_i\right)^{1/2}\right\|_N.}$$

We prove in section~2 that the complex interpolation space
$(A_0,A_1)_\theta$ is the space $N^n$ equipped  with the norm
$$\|(x_1,\ldots, x_n)\|_\theta = \left\|\sum
L_{x_i}R_{x^*_i}\right\|^{1/2}_{B(L_p(\varphi))}$$ where we have
denoted by $L_x$ (resp. $R_x$) the operator of left (resp.  right)
multiplication by $x$ on $L_p(\varphi)$, and where $p={\theta}^{-1}$.
Note that the case $\theta=0$ corresponds to $L_\infty(\varphi)$
identified with $N$ and $\theta=1$ corresponds to $L_1(\varphi)$
identified with $N_*$ in the usual way. Again in the particular case
$N=B(H)$ this result was proved in [P1].  \medskip

We refer to [Ta1] for background on von~Neumann algebras and to [Pa]
for complete boundedness.

We will use several times the following elementary fact.

\proclaim Lemma 0.1. Let $M\subset N$ be von~Neumann algebras. Let
$(p_i)_{i\in I}$ be a directed net of projections in $M$ such that,
for all $ x$ in $M$,  $p_ixp_i$ tends to $x$ in the $\sigma(M,M_*)$
topology. Assume that for each $i$ there is a norm one projection
$P_i\colon \ N\to p_iMp_i$. Then there is a norm one projection $P$
from $N$ onto $M$.

\pf Let ${\cal U}$ be a nontrivial ultrafilter refining the net. For
any $x$ in $N$, we define $P(x) = \lim\limits_{\cal U} P_i(p_ixp_i)$
where the limit is in the $\sigma(M,M_*)$ sense. Then $P(x)\in M$ and
$\|P(x)\|\le \|x\|$. Moreover, for any $x$ in $M$ we have
$$P_i(p_ixp_i)=p_ixp_i.$$ Hence $P(x)=x$ for all $x$ in $M$, and we
conclude that $P$ is a projection from $N$ to $M$.\qed \vfill\eject

\n {\bf \S 1. Projections.}

The main result of this section is the following.

\proclaim Theorem 1.1. Let $M\subset N\subset B(H)$ be von~Neumann
algebras with $M$ semi-finite. If there is a completely bounded ($c.b.$
in short) projection $P\colon \ N\to M$, then there is a norm one
projection $\widetilde P\colon \ N\to M$.

Actually, we use less than complete boundedness, we only need to assume
that there is a constant $C$ such that for all $x_1,\ldots, x_n$ in $N$
we have $$ {\left\|\sum P(x_i)^*P(x_i)\right\| \le C^2\left\|\sum
x^*_ix_i\right\|\quad\hbox{and}\quad \left\|\sum P(x_i)P(x_i)^*\right\|
\le C^2\left\|\sum x_ix^*_i\right\|.}\leqno(1.1)$$ The proof is given
at the end of this section.

\n Notation:\ Let $\varphi$ be a normal faithful semi-finite trace on a
von~Neumann algebra $N$. We denote by $L_2(\varphi)$ the usual
associated Hilbert space. For any $a$ in $N$, we denote by $L_a$ (resp.
$R_a$) the operator of left (resp. right) multiplication by a in
$L_2(\varphi)$, i.e. we set for all $x$ in $L_2(\varphi)$ $$L_ax =
ax,\qquad R_ax=xa.$$ The key lemma in the proof of Theorem~1.1 is the
next statement.

\proclaim Lemma 1.2. Let $N$ be a semi-finite von~Neumann algebra with
a normal faithful semi-finite trace $\varphi$ as above. Consider a
finite set $x_1,\ldots, x_n$ in $N$ and assume
$$\left\|\sum\nolimits^n_1 L_{x_i} R_{x^*_i}\right\|_{B(L_2(\varphi))}
\le 1,\leqno (1.2)$$ then there is a decomposition $x_i=a_i+b_i$ with
$a_i \in N$, $b_i\in N$ such that $$\left\|\left(\sum
a^*_ia_i\right)^{1/2}\right\| + \left\|\sum b_ib^*_i\right\|^{1/2} \le
1.\leqno (1.3)$$

 More generally, the main idea of this paper seems to be the
identification of the expression $$\|(x_1,\ldots, x_n)\| =
\left\|\sum\nolimits^n_1 L_{x_i}R_{x^*_i}\right\|^{1/2}
_{B(L_2(\varphi))}$$ with the norm of a simple interpolation space
obtained by the complex interpolation method. See section 2 for further
details.

\proclaim Corollary 1.3. Let $N$ be as in Lemma~1.2 and let $M$ be a
finite von~Neumann algebra equipped with a normalized finite trace
$\tau$. Let $P\colon \ N\to M$ be any linear map satisfying (1.1). Then
for all finite sequences $x_1,\ldots, x_n$ in $N$ we have
$$\sum^n_1\tau (P(x_i)P(x_i)^*)=\sum^n_1\tau (P(x_i)^*P(x_i)) \le
C^2\left\|\sum^n_1 L_{x_i}R_{x^*_i}\right\|_{B(L_2(\varphi))}.$$

\pf Assume $\big\| \sum L_{x_i}R_{x^*_i}\big\| \le 1$. Let $a_i,b_i$ be
as in Lemma~1.2.  Let us denote $\| x\|_2=(\tau(x^*x))^{1/2}$ for all
$x$ in $M$.  Then we have $$\eqalign{\left(\sum
\|P(x_i)\|^2_2\right)^{1/2} &\le \left(\sum\|P(a_i)\|^2_2\right)^{1/2}
+ \left(\sum\|P(b_i)\|^2_2\right)^{1/2}\cr &\le \left\|\sum P(a_i)^*
P(a_i)\right\|^{1/2} + \left\|\sum P(b_i) P(b_i)^*\right\|^{1/2}\cr
&\le C.}$$\qed

\proclaim Lemma 1.4. Let $N$ be as in  Lemma~1.2 and let $M\subset N$
be a finite von~Neumann subalgebra. Assume that there is a projection
$P\colon \ N\to M$ satisfying (1.1). Then for all nonzero projection
$p$ in the center of $M$ and for all unitary operators $u_1,\ldots,
u_n$ in $M$ we have $$n = \left\|\sum^n_1
L_{pu_i}R_{(pu_i)^*}\right\|_{B(L_2(\varphi))}.\leqno (1.4)$$

\pf Fix $p$ as in Lemma 1.4. By [Ta1, p. 311 ] there is a finite trace
$\tau$ on $M$ with $\tau(p)\not= 0$.  By Corollary 1.3 applied to the
normalized trace $x\to \tau(p)^{-1} \tau(x)$ on $pMp = pM$ we have $$n
= \sum\|pu_i\|^2_2 \le C^2 \left\|\sum^n_1 L_{pu_i}R_{(pu_i)^*}
\right\|_{B(L_2(\varphi))}.$$ To replace $C^2$ by 1 in this inequality,
we use the same trick as Haagerup in [H1]. Let $$T_n  = \sum^n_1
L_{pu_i}R_{(pu_i)^*}.$$ We have for each $k$ $$T^k_n = \sum_{1\le m\le
n^k} L_{x_m}R_{x^*_m}$$ where each $x_m$ is of the form $pu$ with $u$
unitary in $M$. It follows that $$n^k \le C^2\|T^k_n\| \le
C^2\|T_n\|^k$$ hence $n \le C^{2/k}\|T_n\|$. Letting $k$ tend to
infinity we obtain (1.4) (since the other direction is trivial by the
triangle inequality.)\qed

\n {\bf Proof of Lemma 1.2.} We will use the duality between $N^n$ and
$N^n_*$.

\n Let $C$ be the set of elements $(x_i)_{i\le n}$ in $N^n$ which admit
a decomposition $x_i=a_i+b_i$ in $N$ satisfying (1.3). We will show
that if (1.2) holds, then necessarily $(x_i)$ lies in the bipolar
$C^{oo}$ of $C$ in the duality between $N^{n}$ and $N^n_*$. This is
enough to conclude.  Indeed since the set $C$ is clearly convex and
$\sigma(N^n, N^n_*)$ closed we have $C=C^{oo}$, so we obtain that
$(x_i)$ is in $C$ if $(x_i)$ satisfies (1.2).

\n Hence assume given $(x_i)$ satisfying (1.2). Consider $(\xi_i)_{i\le
n}$ in $N^n_*$ and assume\break  $(\xi_i)_{i\le n}\in C^o$. This means
that for any $a_i$ in $N$ such that $${\rm either}\quad
\left\|\sum^n_1 a_ia^*_i\right\|^{1/2}\le 1 \quad {\rm or}\quad
  \left\|\sum a^*_ia_i\right\|^{1/2} \le 1,\leqno (1.5)$$ we have
$$\left|\sum \xi_i(a_i)\right|\le 1.$$ We use the classical
identification $N_* = L_1(\varphi)$ and we use the density of $N\cap
L_1(\varphi)$ in $L_1(\varphi)$. By these well known properties of
$L_1(\varphi)$ for each $\vp>0$ we can find a projection $p$ in $N$
with $\varphi(p)<\infty$ and elements $b_1,\ldots, b_n$ in $pNp$ such
that $$\|\xi_i-b_i\|_{N_*}<\vp.\leqno (1.6)$$ It follows that for any
$(a_i)$ satisfying (1.5) we have $\big|\sum \langle b_i,a_i\rangle\big|
\le 1+n\vp$. So that replacing $b_i$ by ${b_i\over 1+n\vp}$ we may as
well assume (since $\vp>0$ is arbitrary) that, for any $(a_i)$
satisfying (1.5) we have $$\left|\sum \varphi(b_ia_i)\right|\le
1.\leqno (1.7)$$

\n We first claim that this implies $$\varphi\left(\left(\sum
b^*_ib_i\right)^{1/2}\right) \le 1\quad \hbox{and}\quad
\varphi\left(\left(\sum b_ib^*_i\right)^{1/2}\right) \le 1.\leqno
(1.8)$$ Indeed, let $r$ (resp. $c$) be the element of $M_n(N)_*$
corresponding to the $n\times n$ matrix which has coefficients equal to
$b_1,\ldots, b_n$ on the first row (resp. column) and zero elsewhere.
Then by (1.7) $r$ and $c$ are in the unit ball of $M_n(N)_*$. From this
(1.8) immediately follows by the identification between $M_n(N)_*$ and
$L_1(\widetilde\varphi)$ where $\widetilde\varphi$ is the semi-finite
trace defined on $M_n(N)$ by $$\widetilde\varphi((a_{ij})) = \sum
\varphi(a_{ii}).$$ Secondly, we claim that, for any $\delta >0$, $b_i$
can be written as  $b_i=\alpha y_i\beta$ with $\alpha, y_i,\beta$ in
$pNp$ such that $$\varphi(|a|^4) \le 1+\delta\varphi(p), \qquad
\varphi(|\beta|^4) \le 1+\delta\varphi(p)\quad{\rm and }\quad \sum
\varphi(|y_i|^2) \le 1.$$

\n Let $$\alpha = \left(\left(\sum b_ib^*_i\right)^{1/2} +\delta
p\right)^{1/4}\quad \hbox{and}\quad \beta = \left(\left(\sum
b^*_ib_i\right)^{1/2} +\delta p\right)^{1/4}  .$$ Note that we clearly
have $$\alpha^{-2}\left(\sum b_ib_i^*\right) \alpha^{-2} \le \left(\sum
b_ib_i^*\right)^{1/2} \quad \hbox{and}\quad \beta^{-2} \left(\sum
b_i^*b_i\right)\beta^{-2} \le \left(\sum b_i^*b_i\right)^{1/2}.\leqno
(1.9)$$ We also note that $$\varphi(\beta^4) \le
1+\delta\varphi(p)\quad \hbox{and}\quad \varphi(\alpha^4) \le 1+\delta
\varphi(p).\leqno (1.10)$$ Now in the von~Neumann algebra $pNp$ (with
unit $p$) we introduce the analytic $pNp$ valued functions $f_k$
($k=1,...,n$) defined on the strip $S = \{z \in {\bf C}\mid 0<
{Re}(z)<1\}$ by $$f_k(z) = \alpha^{-2(1-z)} b_k\beta^{-2z}.$$ We have
$$f_k(it) =\alpha^{2it}\alpha^{-2} b_k\beta^{-2it}\quad \hbox{and}\quad
f_k(1+it) =\alpha^{2it}b_k\beta^{-2-2it}.$$ Since $\alpha^{2it}$ and
$\beta^{-2it}$ are unitary in $pNp$, it follows that for all real $t$
$$\eqalignno{\sum\|f_k(it)\|^2_{L_2(\varphi)} &=
\varphi\left(\alpha^{-2} \sum b_k b^*_k \alpha^{-2}\right),\cr
\noalign{\hbox{hence by (1.9)  and (1.8)}} &\le \varphi\left(\left(\sum
b_kb^*_k\right)^{1/2}\right) \le 1.}$$
 Note that $f_k$ is bounded on $\overline S$ since $\alpha,\beta$ are
bounded below in $pNp$.  We now invoke    the three lines lemma (cf.
[BL, p. 4]). Note that, as is well known,
 this lemma remains valid for bounded analytic functions on $S$,
 not necessarily continuous on $\overline S$, using the nontangential
boundary  values to extend the functions to $\overline S$.  Using this,
we conclude that, for all $z$ in the strip $S$, we have $$\sum
\|f_k(z)\|^2_{L_2(\varphi)} \le 1.$$ In particular this holds for
$z=1/2$ and we can define $y_k = f_k(1/2)$.  Then we have $b_k = \alpha
y_k\beta$ and all the announced properties hold.  We now return to our
original $n$-tuple $x_1,\ldots, x_n$ in $N$.

\n We have by (1.6) $$\eqalignno{\left|\sum
\langle\xi_k,x_k\rangle\right| &\le \left|\sum \langle
b_k,x_k\rangle\right|+n\vp\cr &\le \left|\sum \varphi(\alpha y_k \beta
x_k)\right| +n\vp\cr \noalign{\hbox{hence by Cauchy-Schwarz and by
(1.10)}} &\le \left(\sum \|\beta x_k
\alpha\|^2_{L_2(\varphi)}\right)^{1/2}+n\vp\cr &=
\left(\varphi\left(\sum \beta
x_k\alpha\alpha^*x^*_k\beta^*\right)\right)^{1/2} + n\vp\cr &=
\left\langle \beta^*\beta, \sum L_{x_k} \alpha\alpha^*
R_{x^*_k}\right\rangle_{L_2(\varphi)}^{1/2} +n\vp\cr &\le \left\|\sum
L_{x_k}R_{x^*_k}\right\|^{1/2}_{B(L_2(\varphi))}
(1+\delta\varphi(p))^{1/2} +n\vp.}$$ Since $\vp,\delta>0$ are
arbitrary, we conclude that if (1.2) holds we have $\big|\sum\langle
\xi_k,x_k\rangle\big| \le 1$ for all $(\xi_k)$ in $C^o$.  Hence we have
$(x_k) \in C^{oo}$ and the proof is complete.  \qed

To prove Theorem~1.1, we will combine Lemma~1.2 with a rather
straightforward extension of some results of Haagerup in [H1] on
injective von~Neumann algebras. Haagerup's work is based on Connes'
ideas on injective factors [Co].

\proclaim Definition 1.5. Let $M\subset N$ be von~Neumann algebras. We
will say that a state $\omega$ on $N$ is an $M$-hypertrace on $N$ if we
have $$\forall\ a\in M\qquad \forall x\in N\quad \omega(ax) =
\omega(xa).$$

\proclaim Theorem 1.6. Let $M\subset N$ be von~Neumann algebras with
$N$ semi-finite and $\varphi$ a faithful normal semi-finite trace on
$N$. The following are equivalent\medskip \item{(i)} $M$ is finite and
there is a norm one projection $P$ from $N$ onto $M$.  \item{(ii)} For
any finite set $u_1,\ldots, u_n$ of unitaries in $M$ and any nonzero
central projection $p$ in  $M$ we have

$$n = \left\|\sum^n_1 L_{pu_i} R_{(pu_i)^*}\right\|_{B(L_2(\varphi))}.
\leqno (1.11)$$

{\sl \item{(iii)}  For every nonzero central projection $p$ in $M$
there is an $M$-hypertrace $\omega$ on $N$, such that $\omega(1-p)=0$.
\item{(iv)} For every state $\omega_0$ on the center of $M$ there is an
$M$-hypertrace $\omega$ on $N$ extending $\omega_0$.} \medskip

\pf The proof of Lemma 2.2 in [H1] extends word for word. We simply
replace there $B(H)$ by $N$ and we denote by $M$ the subalgebra.\qed

\n {\bf Remark.} For the convenience of the reader, we recall the key
idea which is behind the preceding statement. This is best described in
the case when $M$ is a factor. In that case the implication
(ii)~$\Rightarrow$~(i) is proved as follows:\ using the uniform
convexity of $L_2(\varphi)$ one shows that (ii) implies the existence
of a net $(z_\alpha)$ in the unit sphere of $L_2(\varphi)$ such that
$\|uz_\alpha u^*-z_\alpha\|_{L_2(\varphi)}\to 0$ for all $u$ unitary in
$M$. Then if we define on $N$ $$\omega(x) = \lim_{\cal U} \langle
xz_\alpha, z_\alpha\rangle _{L_2(\varphi)} = \lim_{\cal U}
\varphi(xz^*_\alpha z_\alpha)$$ we find that $\omega$ is an
$M$-hypertrace on $N$. Moreover since $M$ is a factor, $\omega$
restricted to $M$ is the trace of $M$. It is then easy to conclude that
there is a norm one projection $P\colon \ N\to M$ which is built
exactly like a conditional expectation.

\n {\bf Proof of Theorem 1.1.} By a well known crossed product
argument, (cf. [Ta2]) there is a semi-finite algebra $\widetilde N$
with $N\subset \widetilde N$ and a completely contractive projection
$Q\colon \ \widetilde N\to N$. Hence, replacing $N$ by $\widetilde N$
we may assume that $N$ is semi-finite.

\n Let $P$ be a projection satisfying (1.1). We first assume $M$
finite. Then by Lemma~1.4, the second assertion in Theorem~1.6 holds.
Therefore, by (ii)~$\Rightarrow$~(i) in Theorem~1.6 there is a norm one
projection from $N$ onto $M$.

\n Now if $M$ is semi-finite, we can write $M = \overline{\cup
p_iMp_i}$ (weak-$*$ closure) where $p_i$ is an increasing net of finite
projections in $M$ such that $p_ixp_i\to x$ in the
$\sigma(M,M_*)$-sense for all $x$ in $M$.  Clearly $x\to p_iP(x)p_i$ is
a projection from $p_iNp_i$ onto $p_iMp_i$ which satisfies (1.1) hence
by the first part of the proof, there is a norm one projection from
$p_iNp_i$ onto $p_iMp_i$. A fortiori there is  a norm one projection
from $N$ onto $p_iMp_i$,
 hence we conclude by Lemma~0.1 that there is a norm one projection
$\widetilde P$ from $N$ onto $M$.\qed \vfill\eject

\n {\bf \S 2. An interpolation theorem.}

Let $N$ be a semi-finite von Neumann algebra, let $1\le p<\infty$ and
let
 $L_p(\varphi)$ be the classical non-commutative $L_p$-space associated
 to a faithful normal semi-finite trace $\varphi$ on $N$. For the
 construction and the basic properties of $L_p(\varphi)$, the classical
references are
 [D,S,Ku,Sti]. For a more concise and recent exposition, see [N].

The following statement extends a result proved in [P1] in the
particular case $N=B(H)$.

\proclaim Theorem 2.1. Fix an integer $n\ge 1$. Let $A_0$ (resp. $A_1$)
be the space $N^n$ equipped with the norm $$\|(x_1,\ldots, x_n)\|_{A_0}
= \left\|\sum^n_1 x_ix^*_i\right\|^{1/2}
\left(\hbox{resp.}\ \|(x_1,\ldots, x_n)\|_{A_1}= \left\|\sum
x^*_ix_i\right\|^{1/2}\right).$$ Then for $0<\theta < 1$, the complex
interpolation space $(A_0,A_1)_\theta$ is the space $N^n$ equipped with
the norm $$\|(x_1,\ldots,x_n)\|_\theta = \left\|\sum
L_{x_i}R_{x^*_i}\right\|^{1/2} _{B(L_p(\varphi))}$$ where $\theta =
1/p$.

\n {\bf Remark.} Note that Theorem~2.1 implies Lemma~1.2 by a well
known property of the interpolation spaces, namely the (norm one)
inclusion $(A_0,A_1)_\theta \subset A_0+A_1$, (see [BL] for more
details). But actually, the proof of Theorem~2.1 is quite similar to
that of Lemma~1.2, although slightly more technical.

We will use Szeg\"o's classical factorization theorem which says that
under a nonvanishing condition, a positive function $W$ in $L_1({\bf
T})$ can always be written as $W=|F|^2$ $(W=\overline FF$ is more
suggestive in view of the non-commutative case) for some $F$ in $H^2$.
Moreover, this can be done with $F$ ``outer'', so that $z\to 1/F(z)$ is
analytic inside the disc, and if we additionally require $F(0)>0$ then
$F$ is unique. Actually, we will need an extension of this theorem (due
to Devinatz) valid for $B(H)$-valued functions. The following
consequence of Devinatz's theorem will be enough for our purposes (cf.
[D], [He]).

\proclaim Theorem 2.2. Let $H$ be a separable Hilbert space and let
$W\colon \ {\bf T}\to B(H)$ be a function such that, for all $x,y$ in
$H$, the function $t\to \langle W(t)x,y\rangle$ is in $L_1({\bf T})$.
Assume that there is $\delta>0$ such that $W(t) \ge \delta I$ for all
$t$.  Then there is a unique analytic function $F\colon \ D\to B(H)$
such that \item{(i)} For all $x$ in $H$, $z\to F(z)x$ is in $H^2(H)$
and its boundary values satisfy almost everywhere on $\bf T$ $$\langle
W(t)x,y\rangle = \langle F(t)x, F(t)y\rangle,$$ \item{(ii)} $F(0)\ge
0$, \item{(iii)} $z\to F(z)^{-1}$ exists and is bounded analytic on
$D$.\medskip

\n The following corollary was pointed out to me by Uffe~Haagerup
during our collaboration on [HP1] (we ended up not using it in our
paper).

\proclaim Corollary 2.3. Consider a von~Neumann subalgebra $N\subset
B(H)$.  Then in the situation of Theorem~2.2, if $W$ is $N$-valued, $F$
necessarily also is $N$-valued.

\pf Indeed, for any unitary $u$ in the commutant $N'$, the function
$z\to u^*F(z)u$ still satisfies the conclusions of Theorem~2.2, hence
(by uniqueness) we must have $F=u^*Fu$, which implies by the
bicommutant theorem that $$F(z) \in N'' =N.$$\qed

\n {\bf Proof of Theorem 2.1.} By well known results, $N$ can be
written as a direct sum of $\sigma$-finite semi-finite algebras. Hence
we can assume that $N$ is $\sigma$-finite and that $H=L_2(\varphi)$ is
separable.

\n Let $\theta=1/p$. Let us denote $L_\infty(\varphi)=N$. Then it is
well known that we have isometrically $$(L_\infty(\varphi),
L_1(\varphi))_\theta = L_p(\varphi).$$ Clearly if $x_1,\ldots, x_n$ in
$N$ are such that $\|(x_1,\ldots, x_n)\|_0 \le 1$, then we have \break
$\big\|\sum L_{x_i}R_{x^*_i}\big\|_{B(L_\infty(\varphi))} \le 1$.
Similarly, it is easy to check by transposition that if\break
$\|(x_1,\ldots, x_n)\|_1 \le 1$, then $\big\|\sum
L_{x_i}R_{x^*_i}\big\|_{B(L_1(\varphi))} \le 1$. Hence, if
$(x_1,\ldots, x_n)$ is in the unit ball of $(A_0,A_1)_\theta$, we have
necessarily by classical interpolation theory $$\left\|\sum
L_{x_i}R_{x^*_i}\right\|_{B(L_p(\varphi))} \le 1$$ where $1/p=\theta$.

\n This is the easy direction. To prove the converse, we assume that
$$\left\|\sum L_{x_i}R_{x^*_i}\right\|_{B(L_p(\varphi))}\le 1.\leqno
(2.1)$$ We will proceed by duality as in the proof of Lemma~1.2. Let
$B$ denote the open unit ball in the space $(A_{0*}, A_{1*})_\theta$.
Note
 that  $A_{0*}$ (resp. $A_{1*}$) coincides with $N^n_*$ equipped with
the norm
 $\|(\xi_1,\ldots, \xi_n)\|=\varphi[(\sum \xi_i^*\xi_i)^{1/2}]$ (resp.
$\|(\xi_1,\ldots, \xi_n)\|=\varphi[(\sum \xi_i\xi_i^*)^{1/2}]$).  Let
$B^o$ be the polar of $B$ in the duality between $N^n$ and $N^n_*$.
 By a well known duality property of interpolation spaces (cf. [BL,Be])
$B^o$ coincides with the unit ball of $(A_0,A_1)_\theta$. Hence to
conclude it suffices to show that (2.1) implies  $(x_1,\ldots, x_n)\in
B^o$.  Equivalently, to complete the proof it suffices to show that, if
(2.1) holds, then for any $(\xi_1,\ldots, \xi_n)$ in $B$ we have
$\big|\sum \xi_i(x_i)\big|\le 1$.  The rest of the proof is devoted to
the verification of this. By density, if we identify again $N_*$ with
$L_1(\varphi)$ in the usual way, we may assume that $\xi_i$ is of the
form $\xi_i(x) = \varphi(b_ix)$ for some $b_i$ in $qMq$ where $q$  is a
finite projection in $M$, i.e. a projection with $\varphi(q)<\infty$.

\n In that case we have $\xi_i(x_i) = \xi_i(qx_iq)$. Note that (2.1)
remains true if we replace $(x_i)$ by $(qx_iq)$. Therefore, at this
point we may as well replace $N$ by the finite von~Neumann algebra
$qNq$ (with unit $q$) so that we are reduced to the finite case. Hence,
for simplicity, we assume in the rest of the proof that $N$ is finite
with unit $I$ and that $\xi_i$ lies in $N$ viewed as a subspace of
$L_1(\varphi)$ (i.e. that the elements $b_i$ above are in $N$ and
$q=I$). By definition of $(A_{0*}, A_{1*})_\theta$, since $(\xi_i)$ is
in $B$ there are functions $f_i\colon \ \overline S\to L_1(\varphi)$
which are bounded, continuous on $\overline S$ and analytic on $S$ such
that denoting $$\partial_0 = \{z\in {\bf C}\mid {Re} z=0\},\qquad
\partial_1 = \{z\in {\bf C}\mid {Re} z =1\}$$ we have
$\xi_i=f_i(\theta)$ for $i=1,...,n$, with $$\sup_{z\in\partial_0}
\varphi\left[\left( \sum f_i(z)^* f_i(z)\right)^{1/2}\right]<1\quad
\hbox{and}\quad \sup_{z\in\partial_1} \varphi\left[\left(\sum f_i(z)
f_i(z)^*\right)^{1/2}\right]<1.\leqno (2.2)$$ Since $\xi_i$ is in
$N\subset L_1(\varphi)$ and $N^n$ is
 dense in $N^n_*$, we may as well assume by a well known fact (cf.
[St]) that the functions $f_1,\ldots, f_n$ take their values into a
fixed finite dimensional subspace of $N\subset L_1(\varphi)$. We are
then in a position to use Theorem~2.2 and its corollary.

\n Let $\delta>0$ to be specified later. We define functions $W_1$ and
$W_2$ on $\partial S  = \partial_0\cup \partial_1$ by setting
$$\eqalign{\forall\ z\in \partial_1\qquad W_1(z) &= \left(\left(\sum
f_i(z) f_i(z)^*\right)^{1/2} +\delta I\right)^{1/2}\cr \forall\ z\in
\partial_0\qquad W_1(z) &= I\cr \forall\ z\in \partial_1\qquad W_2(z)
&= I\cr \forall\ z\in\partial_0\qquad W_2(z) &= \left(\left(\sum
f_i(z)^* f_i(z)\right)^{1/2} +\delta I\right)^{1/2}.}$$ By (2.2) we can
choose $\delta$ small enough so that $$\sup_{z\in\partial_1}
\varphi(W^2_1)<1\quad \hbox{and}\quad \sup_{z\in \partial_0}
\varphi(W^2_2)<1.\leqno (2.3)$$ By Theorem~2.2 and Corollary 2.3,
using  a conformal mapping from $S$ onto $D$, we find bounded
$N$-valued analytic functions $F$ and $G$ on $S$ with (nontangential)
boundary values satisfying $$FF^* = W^2_1\quad \hbox{and}\quad G^*G =
W^2_2.\leqno (2.4)$$ Moreover, $F^{-1}$ and $G^{-1}$ are analytic and
bounded on $S$. Therefore we can write $$f_i(z) = F(z)g_i(z)G(z)$$
where $$g_i(z) = F(z)^{-1} f_i(z)G(z)^{-1}.\leqno(2.5)$$ We claim that
$$\forall\ z\in S\qquad \sum\|g_i(z)\|^2_{L_2(\varphi)}\le 1.\leqno
(2.6)$$ By the three lines lemma, to verify this it suffices to check
it on the boundary of $S$. (Note that we know a priori that
$\sup\limits_{z\in S} \|g_i(z)\|_{L_2(\varphi)}<\infty$ since
$\|F^{-1}\| <\delta^{-1/2}$ and $\|G^{-1}\|\le \delta^{-1/2}$, hence
$g_i$ is an $H^\infty$ function with values in ${L_2(\varphi)}$, and
its nontangential boundary values still satisfy (2.5) a.e. on the
boundary of $S$.)
 We have $$\eqalignno{\forall\ z\in \partial_1\quad \sum
\|g_i(z)\|^2_{L_2(\varphi)} &= \varphi\left(\sum
g_i(z)g_i(z)^*\right)\cr &= \varphi\left(F(z)^{-1} \sum f_i(z)f_i(z)^*
F(z)^{-1*}\right) = \varphi((FF^*)^{-1}(W^2_1-\vp^2I)^2)\cr
\noalign{\hbox{hence by (2.4) and (2.3)}} &\le \varphi(W^2_1)<1.}$$
Similarly, we find $$\forall\ z\in \partial_0 \qquad
\sum\|g_i(z)\|^2_{L_2(\varphi)} \le \varphi(W^2_2)<1.$$ This proves our
claim (2.6).  Finally, if $\theta=1/p$ we have $$L_{2p}(\varphi) =
(N,L_2(\varphi))_\theta \quad \hbox{and}\quad L_{2p'}(\varphi) =
(L_2(\varphi), N)_\theta.$$ Hence by definition of the latter complex
interpolation spaces, since $\|F(z)\|_N = \|W_1\|_N\le 1$ on
$\partial_0$ and (by (2.3)) $\|F(z)\|_{L_2(\varphi)} <1$ on
$\partial_1$, we have $$\|F(\theta)\|_{L_{2p}(\varphi)} \le 1$$ and
similarly $\|G(\theta)\|_{L_{2p'}(\varphi)}\le 1$. Therefore we can
conclude as in section~1:\ we have $\xi_i = f_i(\theta) = F(\theta)
g_i(\theta) G(\theta)$, hence if $(x_i)$ satisfies (2.1) we have by
(2.6) (and Cauchy-Schwarz) $$\eqalign{\left|\sum \xi_i(x_i)\right| &=
\left|\sum \varphi(F(\theta)g_i(\theta) G(\theta)x_i)\right|\cr &\le
\left(\sum \|G(\theta) x_iF(\theta)\|^2_{L_2(\varphi)}\right)^{1/2}\cr
&\le \left\|\sum x_iF(\theta)
F(\theta)^*x^*_i\right\|^{1/2}_{L_p(\varphi)}\cr &\le \left\|\sum
L_{x_i}R_{x^*_i}\right\|^{1/2}_{B(L_p(\varphi))}\le 1.}$$ Thus we have
verified that (2.1) implies $(x_i) \in B^o$. This concludes the proof
of Theorem~2.1.\qed

\vskip24pt \magnification\magstep1
 \baselineskip = 18pt \def\n{\noindent}

 \overfullrule = 0pt \def\qed{{\hfill{\vrule height7pt width7pt
depth0pt}\par\bigskip}}

 \centerline{\bf References}

 \item{[Be]} J. Bergh. On the relation between the two complex methods
of interpolation. Indiana Univ. Math.  Journal 28 (1979) 775-777.

 \item{[BL]} J. Bergh and J. L\"ofstr\"om. Interpolation spaces. An
 introduction. Springer Verlag, New York. 1976.

 \item{[Ca]} A. Calder\'on. Intermediate spaces and interpolation, the
 complex method. Studia Math. 24 (1964) 113-190.

\item{[Co]} A. Connes. Classification of injective factors, Cases
$II_1,II_\infty,III_\lambda,\lambda\neq 1$. Ann.  Math. 104 (1976)
73-116.

\item{[CS]} E. Christensen and A. Sinclair.
 On von Neumann algebras which are complemented subspaces of $B(H)$.
Preprint. To appear.

\item{[D]} A. Devinatz. The factorization of operator valued analytic
functions. Ann. of Math. 73 (1961) 458-495.

\item{[H1]} U. Haagerup. Injectivity and decomposition of completely
 bounded maps in ``Operator algebras and their connection with Topology
 and Ergodic Theory''. Springer Lecture Notes in Math. 1132 (1985)
170-222.

 \item{[H2]} $\underline{\hskip1.5in}$. $L^p$-spaces associated with an
arbitrary von Neumann algebra.  Alg\`ebres d'op\'erateurs et leurs
applications en physique math\'ematique. (Colloque CNRS, Marseille,
juin 1977) Editions du CNRS, Paris 1979.

 \item{[H3]} $\underline{\hskip1.5in}$. Decomposition of completely
 bounded maps on operator algebras. Unpublished manuscript. Sept.
 1980.

\item{[HP1]} U. Haagerup and G. Pisier. Factorization of analytic
functions
 with values in non-commutative $L^1$-spaces and applications. Canadian
 J.  Math. 41 (1989) 882-906.

 \item{[HP2]} U. Haagerup and G. Pisier.  Bounded linear operators
between
 $C^*$-algebras. Duke Math. J. (1993) To appear.

\item{[He]} H. Helson. Lectures on invariant subspaces.  Academic
Press, New-York 1964.

\item {[KR]} R. Kadison  and J. Ringrose.    Fundamentals of the theory
of operator algebras, Vol. II, Advanced Theory,    Academic Press,
New-York      1986.

\item {[Ku]} R. Kunze. $L_p$ Fourier transforms on locally compact
unimodular groups. Trans. Amer. Math. Soc. 89 (1958) 519-540.

\item {[N]} E. Nelson. Notes on non-commutative integration.  J. Funct.
Anal. 15 (1974) 103-116.

 \item{[Pa]} V. Paulsen.  Completely bounded maps and dilations. Pitman
Research Notes 146.  Pitman Longman (Wiley) 1986.

 \item{[P1]} G. Pisier.  Espace de Hilbert d'op\'erateurs et
interpolation complexe. Comptes Rendus Acad. Sci. Paris 316 (1993)
4-52.

 \item{[P2]} The operator Hilbert space $OH$, complex interpolation and
tensor norms.  Memoirs Amer. Math. Soc. (submitted)

\item{[P3]}  $\underline{\hskip1.5in}$. Remarks on complemented
subspaces of $C^*$-algebras. Proc. Roy. Soc.  Edinburgh, 121A (1992)
1-4.

\item {[P4]} $\underline{\hskip1.5in}$. Complex interpolation and
regular operators between Banach lattices. Arch. der Mat. (Basel) To
appear.

\item{[Ru]} Z.J. Ruan. Subspaces of $C^*$-algebras. J. Funct. Anal. 76
 (1988) 217-230.

\item{[S]} I. Segal. A non-commutative extension of abstract
integration. Ann. of Math. 57 (1953) 401-457.

 \item{[St]} J. Stafney. The spectrum of an operator on an
interpolation space.  Trans. Amer. Math. Soc. 144 (1969) 333-349.

\item{[Sti]} W. Stinespring. Integration theorems for gages and duality
for unimodular groups. Trans. Amer. Math.  Soc. 90 (1959) 15-26.

\item{[Ta1]} M. Takesaki. Theory of Operator Algebras I.
Springer-Verlag New-York 1979.

\item{[Ta2]}$\underline{\hskip1.5in}$. Duality for crossed products
 and the structure of von Neumann algebras of type III.  Acta Math. 131
(1973) 249-308.

\item{[Te1]} M. Terp. Interpolation spaces between a von Neumann
algebra and its predual. J. Operator Th.  8 (1982) 327-360.

\item{[Te2]} $\underline{\hskip1.5in}$. $L^p$-spaces associated with
von Neumann algebras. Preprint. Copenhagen University. June 1981.

\item{[T]} J.Tomiyama. On the projection of norm one in $W^*$-algebras.
Proc. Japan Acad. 33 (1957) 608-612.

\vskip12pt

\vskip12pt

Texas A\&M University

College Station, TX 77843, U. S. A.

and

Universit\'e Paris VI

Equipe d'Analyse, Bo\^\i te 186,

75252 Paris Cedex 05, France

 \end

 Let $A_0$ (resp. $A_1$) be $N^n$ equipped with the norm
$$\eqalign{\|(x_1,\ldots, x_n)\|_{A_0} &= \left\|\sum
x_ix^*_i\right\|^{1/2}\cr \big(resp.\quad \|(x_1,\ldots, x_n)\|_{A_1}
&= \big\|\sum x^*_xx_i\big\|^{1/2}\big).}$$

We will prove in section~2 that the norm in the complex interpolation
space $(A_0,A_1)_{1/2}$ is exactly the norm given by (1.4). We will
also describe $(A_0,A_1)_\theta$ for $0<\theta<1$.  We may as well
assume by density that $b_i\in pNp$ $(i=1,\ldots, n)$ for a projection
$p \in N$ with $\varphi(p)<\infty$.